\documentclass[11pt]{amsart}

\usepackage{amssymb,upgreek}

\usepackage{url}
\usepackage{bm}
\usepackage[left=1in,top=1in,right=1in,bottom=1in,head=.2in]{geometry}
\usepackage{enumerate}

\usepackage[textsize=scriptsize]{todonotes}
\usepackage{color}
\usepackage{mathtools} % for overbraces

\usepackage{fancyhdr}
\usepackage{mathrsfs}
\usepackage[cmtip,all]{xy}
\usepackage{hyperref}

\swapnumbers
\newtheorem{theorem}{Theorem}[section] % numbered theorems, lemmas, etc

\newtheorem{proposition}[theorem]{Proposition}
\newtheorem{corollary}[theorem]{Corollary}
\newtheorem*{theorem*}{Theorem}

\newtheorem*{fcthm*}{Finite Cork Theorem}
\newtheorem*{ccthm*}{Cork Consolidation Theorem}
\newtheorem*{csthm*}{Cork Separation Theorem}
\newtheorem*{icthm*}{Infinite Cork Theorem}
\newtheorem*{icthms*}{Infinite Cork Theorems}
\newtheorem*{aclemma*}{\ac-Lemma}
\newtheorem*{mclemma*}{Multicork Lemma}
\newtheorem*{multicorktheorem*}{Multicork Theorem}
\newtheorem*{lemma*}{Lemma}
\newtheorem*{corollary*}{Corollary}
\newcommand{\thistheoremname}{}
\newtheorem{genericthm}[theorem]{\thistheoremname}

\theoremstyle{definition}
\newtheorem{definition}[theorem]{Definition}
\newtheorem{remark}[theorem]{Remark}

\newtheorem*{remark*}{Remark}
\newtheorem*{definition*}{Definition}
\newtheorem*{remarks*}{Remarks}
\newtheorem*{addenda*}{Addenda}

\newcommand{\fig}[3]{\begin{figure}[h!] \includegraphics[height=#1pt]{#2}#3\end{figure}}

\newcommand{\bit}[1]{\textbf{\textit{#1}}} % {\textit{#1}} %

\newcommand{\bc}{\mathbb C}

\newcommand{\st}{\ \vert\ }

\newcommand{\pt}{\textup{pt}}

\newcommand{\sto}{\!\!\xymatrix@C=1em{{}\ar@{~>}[r]&{}}\!\!}

\newcommand{\del}{\partial}

\newcommand{\sss}{S^2\!\!\times\!S^2}

\newcommand{\cponebar}{\smash{\overline{\bc P}^1}}
\newcommand{\cptwo}{\bc P^2}
\newcommand{\cptwobar}{\smash{\overline{\bc P}^2}}
\newcommand{\interior}{\textup{int}}
\newcommand{\ac}{\textup{AC}}

\newcommand{\overbar}[1]{\mkern 1mu\overline{\mkern-2.5mu#1\mkern-1mu}\mkern 1mu}

\newcommand{\cs}{\mathop\#}

\newcommand{\foot}[1]{\setcounter{footnote}{1}\footnote{\ #1}}

\newcommand{\items}{\begin{itemize}[leftmargin=25pt,rightmargin=5pt]
  \setlength\itemsep{2pt}\vskip-3pt\vskip-3pt}
\newcommand{\stopitems}{\end{itemize}}

\address{Max Planck Insitut f\"ur Mathematik,
Vivatsgasse 7, 53111 Bonn, Germany}
\email{hrschwartz@mpim-bonn.mpg.de} 

\begin{document}

\def\A{\overbar A}
\def\H{\overbar H}

\parskip 2pt

\title{A note on the complexity of \lowercase{h}-cobordisms}
\author{Hannah R. Schwartz}

\maketitle

\begin{abstract}
We show that the number of double points of smoothly immersed $2$-spheres representing certain homology classes of an oriented, smooth, closed, simply-connected $4$-manifold $X$ must increase with the complexity of corresponding h-cobordisms from $X$ to $X$. As an application, we give results restricting the minimal number of double points of immersed spheres in manifolds homeomorphic to rational surfaces. 
\end{abstract}

\section{Introduction} 

It follows from work of Wall \cite{wall:forms}, \cite{wall:forms1}, \cite{wall:4-manifolds}, \cite{wall:diffeomorphisms} that any two oriented, smooth, closed, simply-connected $4$-manifolds that are homotopy equivalent cobound a smooth ($5$-dimensional) h-cobordism. By Freedman \cite{freedman:simply-connected}, such an h-cobordism is topologically a product and so these manifolds are in fact homeomorphic. That h-cobordisms need not smoothly be products is due originally to Donaldson \cite{donaldson}, who gave the first examples of homeomorphic but not diffeomorphic smooth $4$-manifolds, motivating the extensive study of corks,  plugs, and other phenomena relating the distinct smooth structures on a given topological $4$-manifold (see for instance \cite{akbulut:contractible},  \cite{curtis-freedman-hsiang-stong}, \cite{matveyev}). 

In the 90's, Morgan and Szab\'o \cite{morganszabo} produced h-cobordisms \emph{between diffeomorphic $4$-manifolds} that are not products, and in fact, have arbitrarily high complexity\foot{The results in \cite{morganszabo} generalize earlier results, also due to Morgan and Szab\'o, from \cite{morganszabo2}.}(see Definition \ref{complexity1}).  Given a closed, simply-connected $4$-manifold\foot{We work throughout in the smooth, oriented category.} $X$, our main result compares the complexity of h-cobordisms from $X$ to $X$ associated to a class $\sigma \in H_2(X)$ to the minimum number of double points of any immersed sphere representing that class. We refer to this minimum as the complexity $c_\sigma$ of the class $\sigma$ (see Definition \ref{complexity2}), and prove the following.  

\begin{theorem}\label{main1} 
Let $\sigma \in H_2(X)$ be a class of square $\pm 1$ or $\pm 2$. The complexity of the h-cobordism from $X$ to $X$ corresponding to the reflection $\rho_\sigma:Q_X \to Q_X$ is less than or equal to
\begin{enumerate}[(i)]
\item $2c_\sigma+1$ when $\sigma^2= \pm 2$, and
\item $4c_\sigma+1$ when $\sigma^2 = \pm 1$. 
\end{enumerate}
\end{theorem}

Combining Theorem \ref{main1} with Morgan and Szab\'o's results on complexity gives a new method of finding classes with arbitrarily high complexity in smooth manifolds homeomorphic but not necessarily diffeomorphic to a rational surface $\mathbb{C}P^2 \cs n \cptwobar$. This allows us to generalize known results on the complexity of homology classes in $\mathbb{C}P^2 \cs n \cptwobar$ (for instance, in this setting Corollary \ref{main2} is implied by Ruberman \cite{danny}) to all smooth manifolds in this homeomorphism class. 

\begin{corollary}\label{main2}
Let $Y_1, Y_2, Y_3,\dots$ be smooth $4$-manifolds with $Y_n$ homeomorphic to $\mathbb{C}P^2 \cs n \cptwobar$. For each $c \in \mathbb{N}$ and $k\geq-2$, infinitely many $Y_n$ have a characteristic class of square $k$ and complexity at least $c$. 
\end{corollary} 

\begin{corollary}\label{main3}
Let $Y$ be a smooth $4$-manifold homeomorphic to a rational surface. For each $c \in \mathbb{N}$, there is a finite upper bound on the square of characteristic homology classes with complexity $c$. 
\end{corollary}

Both corollaries supplement a large collection of related work, starting in the 80's with Kuga \cite{kuga} and Suciu \cite{suciu} who appealed to Donaldson's results \cite{donaldson} to bound the complexity and minimal genus of classes in $S^2 \times S^2$ and $\mathbb{C}P^2$. A decade later, Donaldson invariants were used to identify genus minimizing surfaces in symplectic $4$-manifolds \cite{km1}, \cite{morganszabotaubes}-- in particular, the Thom conjecture was solved by Kronheimer and Mrowka \cite{km2}. Later contributions discussing complexity and minimal genus in rational surfaces and other connected sums $p\cptwo \cs q\cptwobar$ are quite vast, and include \cite{fs:immersions}, \cite{gg}, \cite{gqz}, \cite{gao},  \cite{lawson}, \cite{LL}, \cite{danny}, and \cite{saso}. 

\smallskip
\noindent
{\bf Acknowledgments.} Thank you to both Paul Melvin and Danny Ruberman for helpful discussions, as well as to the referee for their detailed comments. The author is also grateful to the supportive math communities at both the MPIM and Bryn Mawr College, where this work was carried out. 

\section{Preliminaries} 

Let $X$ be a simply-connected, closed $4$-manifold with a class $\sigma \in H_2(X)$ of square $\pm 1$ or $\pm 2$. 

\begin{definition} \label{reflect} An h-cobordism $W$ from $X$ to $X$ is said to \bit{correspond to} an isometry $\phi$ if the ``lower" inclusion $\alpha: X \hookrightarrow W$ and the ``upper" inclusion (which reverses orientation) $\beta: X \hookrightarrow W$ have induced maps on homology with $\beta_*^{-1}\circ \alpha_*= \phi$. For each $\phi$, there is a unique such h-cobordism from $X$ to $X$ up to diffeomorphism rel boundary. This follows from work of Lawson \cite{lawson} as well as Stallings \cite{stallings}; Kreck later offered a different proof of this fact in \cite{kreck}. We will be mostly concerned with isometries $\rho_\sigma\colon Q_X \to Q_X$ corresponding to classes $\sigma \in H_2(X)$ of square $\pm 1$ or $\pm 2$, given by $\rho_\sigma(\xi) = \xi \mp 2(\xi \cdot \sigma)\sigma$ in the first case and $\rho_\sigma(\xi) =  \xi \mp (\xi \cdot \sigma)\sigma$ in the latter. Such an isometry is called the \bit{reflection in $\sigma$}. 
\end{definition}

\begin{remark*}  If $\sigma$ is represented by an {\it embedded} sphere $\Sigma$, then the associated reflection is induced by an automorphism of $X$; therefore the corresponding h-cobordism is diffeomorphic to a product. Such a diffeomorphism of $X$ may be defined as follows.  On a tubular neighborhood $N_\Sigma$ of $\Sigma$, reflect in $\Sigma$ and in each normal disk. Since $\partial N_\Sigma$ is diffeomorphic to either $S^3$ or $L(2,1)$, by classical results we can insert an isotopy from the reflection to the identity map on a collar of $\partial N_\Sigma$. Extending by the identity on the rest of $X$ then obtains the desired automorphism. 
\end{remark*}

For a class $\sigma$ \emph{not} represented by an embedded sphere, the reflection $\rho_\sigma$ need not be realized by an automorphism of $X$; instead one must construct a more complicated h-cobordism from $X$ to $X$ corresponding to $\rho_\sigma$. 

\begin{remark} \label{construction} The h-cobordism $W$ corresponding to a reflection $\rho_\sigma$ in a class $\sigma$ can be constructed explicitly as follows.   Choose a disk $D$ embedded in $X$.  The normal disk bundle $\nu_{\partial D}$ of the boundary circle $\del D$ acquires a natural framing from the unique framing of the normal disk bundle $\nu_D$ of $D$.  Let $Z$ be the $5$-manifold obtained from $X\times I$ by attaching a $2$-handle along the copy of $\del D$ in $X\times\{1\}$ with this framing.  Then $Z$ is a cobordism from $X$ to a $4$-manifold $X^\circ$ diffeomorphic to $X\cs\sss$ (where the connected sum is performed along the $4$-ball $\nu_D \cong D^2 \times D^2$), with quadratic form $Q_{X^\circ}$ naturally identified with $Q_X\oplus H$.  Here $H$ is the hyperbolic form generated by $a$ and $a^*$ with $a^2=(a^*)^2 = 0$ and $a\cdot a^* = 1$, where $a$ is represented by the disk $D \cap (X -\nu_{\partial D})$ capped off with a parallel copy of the core of the $2$-handle, and $a^*$ is represented by the belt sphere of the $2$-handle.  

Now form $W$ by gluing $Z$ to $-Z$ along $X^\circ$ using an (orientation preserving) automorphism $\phi: X^\circ \to X^\circ$.inducing an isometry $\Phi: Q_X \oplus H \to Q_X \oplus H$ such that
\begin{enumerate}
\item $\Phi(a^*) \cdot a^* = \pm 1$
\item $\pi \circ \Phi$ restricts to $\rho_\sigma$ on $Q_X$
\end{enumerate} 

\noindent where $\pi: Q_X \oplus H \to Q_X$ denotes projection. Condition (1) ensures that the ascending sphere of the 2-handle and attaching sphere of the 3-handle pair algebraically\foot{Or equivalently, $H_*(W, X)=0$ for each component $X \subset \partial W$.} (and so $W$ is an h-cobordism), while Condition (2) guarantees that $W$ corresponds to $\rho_\sigma$. 
\end{remark}

It follows from Smale \cite{smale} that, like the h-cobordisms constructed in Remark \ref{construction} above, \emph{any} h-cobordism $W$ from $X$ to $X$ (and more generally, between any pair of closed, simply-connected $4$-manifolds) can be built as a handlebody from $X \times I$ using only handles of index $2$ and $3$. In other words, $W$ may be obtained by first constructing a cobordism from $X$ to $X \cs n S^2 \times S^2$ by attaching $n$ $2$-handles along embedded curves in $X \times \{1\} \subset X \times I$, and then attaching $n$ $3$-handles along embedded $2$-spheres in $X \cs n S^2 \times S^2$. There are then two sets of distinguished 2-spheres smoothly embedded in $X \cs n S^2 \times S^2$, namely the ascending spheres $A_1, \dots, A_n$ of the 2-handles, and the attaching spheres $B_1, \dots, B_n$ of the 3-handles. Since $W$ is an h-cobordism, these handles can be slid over each other until $A_i \cdot B_j = \delta_{ij}$ (and perturbed if necessary so that each pair of spheres intersects transversally). 

\begin{definition}\label{complexity1} 
Consider the sum $$\sum_{i=1}^n \sum_{j=1}^n  |A_i \pitchfork B_j|- \delta_{ij},$$ i.e. the number of ``excess" intersection points between the spheres which algebraically cancel. The minimum value of this sum over all handlebody structures for $W$ with only $2$ and $3$-handles is called the \bit{complexity of the h-cobordism $W$}, or of its corresponding isometry (see Definition \ref{reflect}). 
\end{definition}

Observe that an h-cobordism has complexity zero if and only if it can be built from $X \times I$ without adding any handles (and so is diffeomorphic to a product). However, the complexity may be high even for h-cobordisms admitting handlebody structures with only a single $2$ and $3$-handle. Whether there exist smooth, closed, simply-connected $4$-manifolds for which \emph{every} h-cobordism between them requires more than one $2$ and $3$-handle is currently unknown. The complexity of h-cobordisms was initially studied by Morgan and Szab\'o in \cite{morganszabo}, who use the Seiberg-Witten invariants of $4$-manifolds with $b_2^+=1$ to produce h-cobordisms of arbitrarily high complexity\foot{The assumption on the second Betti number is critical to their arguments, which rely on the chamber structure of the Seiburg-Witten invariant in this case.}. We relate their observations to the study of the complexity of surfaces.

\begin{definition}\label{complexity2} 
Every self-transverse, smoothly immersed sphere $\Sigma$ in an oriented $4$-manifold $X$ has a \bit{self-intersection} number $\mathcal I(\Sigma) \in \mathbb{Z}[\pi_1(X)]$ originally defined by Whitney in his $1944$ paper \cite{whitney}. In the case when the ambient manifold $X$ is simply-connected, $\mathcal I(\Sigma) \in \mathbb{Z}$ is simply the signed sum of the double points of the immersion. The \bit{complexity $c_\sigma$ of a class $\sigma \in H_2(X)$} is the minimum self-intersection taken over all smoothly immersed spheres representing that class (this minimum is well-defined, by the Hurewicz theorem). 
\end{definition}

In the next section, we will explicitly produce h-cobordisms whose complexity depends on the complexity of immersed spheres used  in their construction. It will be necessary to understand the following automorphisms of $4$-manifolds examined by Wall in \cite{wall:diffeomorphisms}. 

\begin{definition} \label{maps} Consider the spheres $\mathcal A = S^2 \times \pt$ and $\mathcal A^* = \pt \times S^2$ in $X^\circ$, now considered literally equal to $X\cs\sss$. Let $\Sigma$ be a self-transverse, smoothly immersed $2$-sphere in $X^\circ$ representing a class $\sigma \in H_2(X^\circ)$. The following \bit{elementary automorphisms} $X^\circ \to X^\circ$ will be critical to our arguments. Although these maps are defined using the immersion $\Sigma$, we index them by the homology class $\sigma$ since we are concerned only with their induced maps on homology (or equivalently, as in Definition \ref{reflect}, the diffeomorphism type of the corresponding h-cobordism). Indeed, up to isotopy these maps may depend on the choice of the immersion $\Sigma$, and in some cases, an arc in $X^\circ - \Sigma$. 

\noindent {\bf (a)} \emph{When $\sigma^2=2 \ell$ and $\Sigma$ is immersed in $X^\circ - (\mathcal A \cup \mathcal A^*)$, the maps $E_\sigma$ and $E^*_\sigma$:} We define $E_\sigma$ explicitly, and $E^*_\sigma$ analogously with the roles of $\mathcal A$ and $\mathcal A^*$ interchanged. First surger $\mathcal A^*$, replacing its neighborhood by a neighborhood $N_\gamma$ of an embedded circle $\gamma \subset X$ bounding the hemisphere $D$ of the sphere $\mathcal A$ that remains in $X$. Note that the surgery uniquely determines (up to isotopy) one of two possible trivializations $f: S^1 \times B^3 \to N_\gamma$ corresponding to an element of $\pi_1(SO(3)) \cong \mathbb{Z}_2$.

Now, take two copies of $D$, and tube their interiors to distinct points in $\Sigma$ along disjointly embedded arcs in $X -\Sigma$. This gives an immersed annulus $A: S^1 \times I \looparrowright X$, thought of here as a homotopy of $\gamma$. After reparametrizing if necessary, $A$ defines an isotopy of $\gamma$ which extends to an ambient isotopy $\psi: X \times I \to X$ that ``sweeps" the curve $\gamma$ across the sphere $\Sigma$ and then back to its original position. The end of this ambient isotopy gives an automorphism $\psi_1: X \to X$ fixing $N_\gamma$ setwise, and hence a new trivialization $f'= \psi_1 \circ f : S^1 \times B^3 \to N_\gamma$ of the neighborhood $N_\gamma$. 

To see which framing of $\gamma$ is given by $f'$, consider the sub-bundle $N^\circ_\gamma \subset N_\gamma$ with fibers the 2-disks normal to the annulus $A$, also preserved setwise by $\psi_1$. Since $\sigma^2= 2\ell$, the restriction $\psi_1|_{N^\circ_\gamma}$ is a bundle map corresponding to an even element of $\pi_1(SO(2)) \cong \mathbb{Z}$. This implies that the bundle map $\psi_1 |_{N_\gamma}$ corresponds to the trivial element of $\pi_1(SO(3))$, i.e. the trivializations $f$ and $f'$ are isotopic. Thus (after an isotopy in a collar of $\partial N_\gamma$), $\psi_1$ may be assumed to fix $N_\gamma$ not only setwise but \emph{pointwise}, allowing $\psi_1|_{X^\circ - N_\gamma}$ to extend to an automorphism $E_\sigma$ of $X^\circ$. The effect of $E_\sigma$ on the intersection form of $X^\circ$ is analyzed by Wall in \cite[Corollary 2, Section 3]{wall:diffeomorphisms}. In particular, the diffeomorphism $E_\sigma$ induces the isometry $Q_{X^\circ} \to Q_{X^\circ}$ sending $$a \mapsto a + \sigma - \ell a^*, ~~~a^* \mapsto a^*, ~~~\xi \mapsto \xi - (\xi \cdot \sigma) a^*$$
for each $\xi \in H_2(X)$. 

\noindent {\bf (b)} \emph{The map $R$:} On a tubular neighborhood of $\mathcal A \cup \mathcal A^*$, define $R$ as the reflection across the diagonal of $S^2 \times S^2$. This map exchanges $\mathcal A$ and $\mathcal A^*$, and restricts to a map on the boundary $S^3$ that may be isotoped by Cerf \cite{cerf} to the identity and then extended across the rest of $X^{\circ}$ by the identity. The induced automorphism of $Q_{X^\circ}$ sends $$a \mapsto a^*, ~~~a^* \mapsto a, ~~~\xi \mapsto \xi$$
for each $\xi \in H_2(X)$.

\noindent {\bf (c)} \emph{When $\sigma^2 = \pm1$ and $\Sigma$ is embedded, the map $S_\sigma$:} On a tubular neighborhood of $\Sigma$ diffeomorphic to a punctured copy of $\pm \mathbb{C}P^2$, define $\mathcal S_\sigma$ as complex conjugation; this restricts to a map on the boundary $S^3$ that may be isotoped by Cerf \cite{cerf} to the identity and then extended across the rest of $X^{\circ}$ by the identity. The induced automorphism of $Q_{X^\circ}$ is the reflection $\rho_\sigma$ and so each class $\xi \in H_2(X^\circ)$ is sent to $\xi-2(\xi \cdot \sigma) \sigma$.   

\end{definition}

\begin{remark}\label{Esigma} 
Since the automorphism $\psi_1$ of $X$ constructed in part (a) of the definition above fixes $N_\gamma$ pointwise, the surgery to obtain $X^\circ$ may be performed on a slightly smaller neighborhood of $\gamma$ so that the elementary automorphism $E_\sigma$ restricts to the identity on a neighborhood of the sphere $\mathcal A^*$. 
\end{remark}

\begin{remark}\label{Esigma2}
Suppose that the self-intersection $\mathcal I(\Sigma)$ of the immersion $\Sigma$, as in Definition \ref{complexity2}, is equal to $-\sigma^2/2= -\ell$. In this special case, one can understand the map on homology induced by $E_\sigma$ geometrically by following the disk $D$ during the isotopy $\psi$ as it is dragged along with its boundary $\gamma$ around $X$. To avoid intersecting $\gamma$ as it crosses a double point of the immersion $A$, the interior of $D$ is forced to wrap around the boundary of a $3$-ball normal to $\gamma$, tubing one copy of $\mathcal A^*$ to $D$ for each double point of $\Sigma$, as illustrated in Figure \ref{diffeo}. The fact that $\mathcal I(\Sigma) = -\ell$ implies that the intersection of the surface $\Sigma$ with a generic push-off vanishes\foot{This follows from the well-known formula $2\mathcal I(\Sigma)= \Sigma \cdot \Sigma - \sigma^2$ due to Lashof and Smale \cite{lashofsmale}.}; therefore the bundle map $\psi_1|_{N^\circ_\gamma}$ corresponds to the trivial element of $\pi_1(SO(2))$. Hence, the final isotopy of $\psi_1$ (supported in a collar of $\partial N_\gamma$ and arranging $N_\gamma$ to be pointwise fixed) can be done preserving the sub-bundle $N^\circ_\gamma$ so that no additional copies of the boundary of $3$-ball normal to $\gamma$ are tubed to the disk $D$.   
\end{remark}

%%%%%%%%%% FIG 1 %%%%%%%%%%
\fig{300}{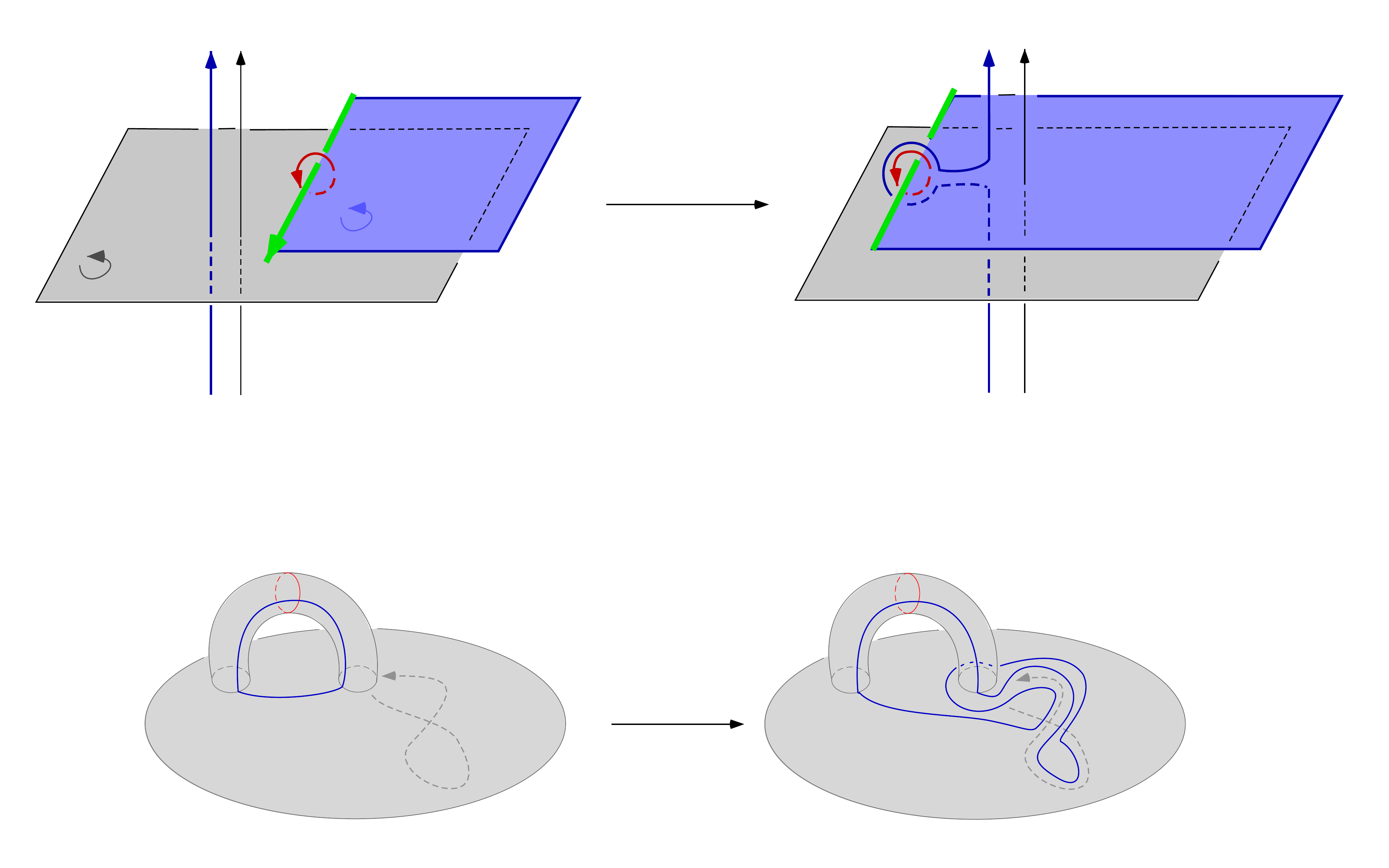}{
\put(-410,130){The 3D analogy:}
\put(-360,271){$\gamma$}
\put(-335,235){$D$}
\put(-153,273){$\gamma$}
\put(-394,180){$A$}
\put(-430,220){$A$}
\caption{Dragging the curve $\gamma$ and the disk $D$ across $\Sigma$ using the isotopy defined by $A$. The red circle represents the sphere $\mathcal A^*$, which bounds a 3-ball normal to $\gamma$.}
\label{diffeo}}
%%%%%%%%%%%%%%%%%%%%%%%%%%%

\section{Results}

Consider the spheres $\mathcal A$ and $\mathcal A^*$ representing classes $a,a^* \in H_2(X^\circ)$ as in the previous section.

\smallskip

\noindent \emph{Proof of Theorem \ref{main1}:}
 Let $\Sigma \subset X$ be a self-transverse, smoothly immersed sphere representing the class $\sigma$. Assume that the immersion $\Sigma$ has complexity $c_\sigma$ and is disjoint from the disk $D \subset X$, so that it may also be thought of in $X^\circ$. We argue the cases $(1)$ and $(2)$ separately.
\smallskip

\noindent
\bit{Case 1: $\sigma^2= \pm 2$.}
\smallskip

Begin by adding ``kinks" as in Figure \ref{kink} to the immersion $\Sigma$, until the self-intersection $\mathcal I(\Sigma)= \mp 1$ (see Remark \ref{Esigma2}). Since at most $c_\sigma +1$ extra kinks are needed to achieve this, the immersion $\Sigma$ now has $n \leq 2c_\sigma +1$ double points. Using this modified immersion, construct a cobordism $W$ from $X$ to $X$ as in Remark \ref{construction}, with $\phi=E_\sigma E^*_{\pm \sigma} E_\sigma$ the composition of the elementary automorphisms of $X^\circ$ from  part (a) of Definition \ref{maps}. We claim that $W$ is the h-cobordism corresponding to the reflection $\rho_\sigma$, and verify this by checking that $\phi$ satisfies both Condition (1) and (2) from Remark \ref{construction}. First, note that the induced map on $H_2(X^\circ)$ sends $a^* \mapsto \mp a$. Furthermore, setting $a,a^*=0$, the induced map sends $\xi \mapsto \xi \mp (\xi \cdot\sigma)\sigma$ for each $\xi \in H_2(X)$, which gives $\rho_\sigma$ as desired\foot{Basically the same fact was noticed by Wall in \cite{wall:forms1}, and earlier by Hasse \cite{hasse}.}

%%%%%%%%%% FIG 1 %%%%%%%%%%
\fig{90}{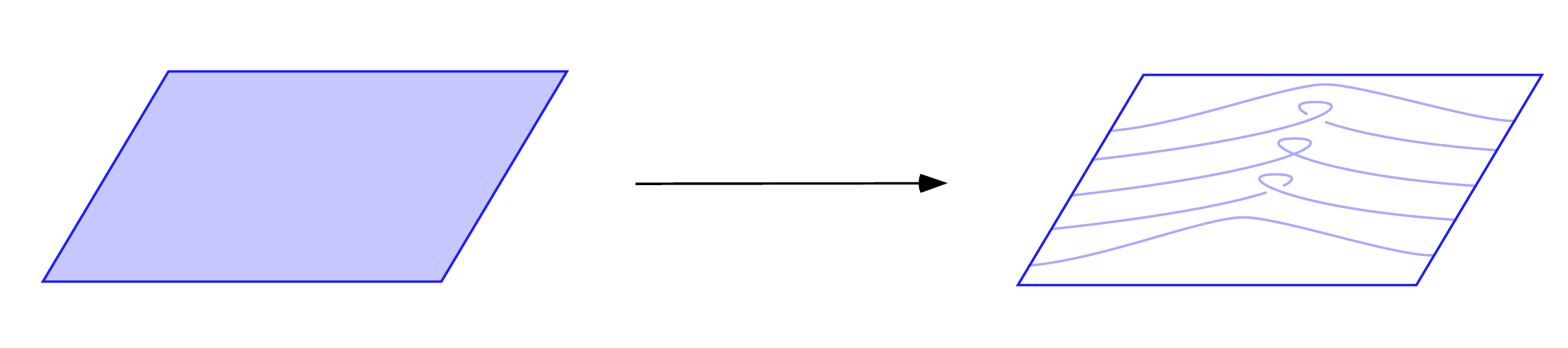}{
\put(-230,50){Add a kink}
\caption{Creating a kink in an immersion adds a single double point, which can be made to have either sign.}
\label{kink}}
%%%%%%%%%%%%%%%%%%%%%%%%%%%

Recall from Definition \ref{complexity1} that the complexity of $W$ is bounded above by the number of ``excess" intersection points between the spheres $\phi(\mathcal A^*)$ and $\mathcal A^*$  in $X^\circ$ (the ascending sphere of the $2$-handle of $W$ and the attaching sphere of the $3$-handle, respectively). To compute this geometric intersection number, we follow the image of the sphere $\mathcal A^*$ under each elementary automorphism of $X^\circ$ in the composition $\phi$. 

\begin{enumerate}
\item By Remark \ref{Esigma}, the map $E_\sigma$ fixes $\mathcal A^*$. 

\item The map $E^*_{\pm \sigma}$ sends $\mathcal A^*= E_\sigma (\mathcal A^*)$  to a sphere $\mathcal A^*_1$ gotten from the connected sum\foot{Recall from Definition \ref{maps} that each elementary automorphism $E_\sigma$ depends on a choice of paths along which to tube copies of the disk $D$ to the immersion $\Sigma$. Different tube choices do not affect the h-cobordism $W$ up to diffeomorphism (as remarked in Definitions \ref{reflect} and \ref{maps}), but do determine different embedded arcs in the complement of $\mathcal A^* \cup \Sigma$ along which the connected sum $\mathcal A^* \cs \pm \Sigma$ is formed in Step (2) above.} $\mathcal A^* \cs \pm \Sigma$ by using the Norman trick \cite{norman} repeatedly to eliminate its $n$ double points by ``tubing" them to parallel copies of $\mathcal A$, as in Figure \ref{spheres}. This produces an embedded sphere $\mathcal A^*_1$ intersecting $\mathcal A^*$ transversally in $n$ points, one for each parallel copy of $\mathcal A$ used in the construction. 

\item Again by Remark \ref{Esigma}, the map $E_\sigma$ restricts to the identity on a neighborhood of $\mathcal A^*$.  Therefore, the sphere  $\mathcal A^*_1$ and its image $\mathcal A^*_2= E_\sigma(\mathcal A^*_1)$ intersect $\mathcal A^*$ in the same number of points, namely $n$. 

\end{enumerate}

%%%%%%%%%% FIG 1 %%%%%%%%%%
\fig{180}{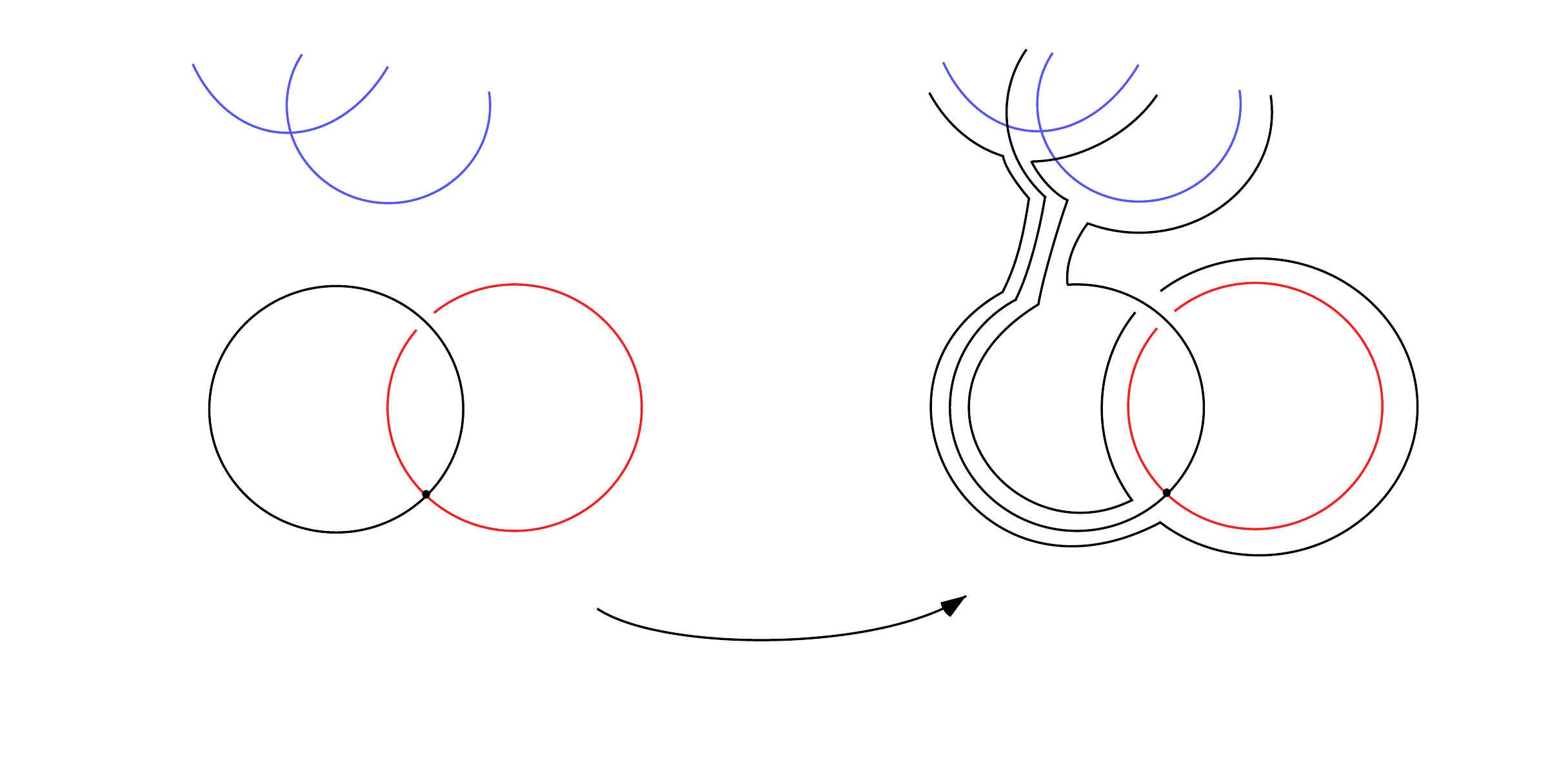}{
\put(-254,137){$\Sigma$}
\put(-214,90){$\mathcal A$}
\put(-335,90){$\mathcal A^*$}
\put(-170,90){$\mathcal A_1^*$}
\put(-207,12){Step (2)}
\caption{A schematic of Step (2) in the proof of Theorem \ref{main1}.}
\label{spheres}}
%%%%%%%%%%%%%%%%%%%%%%%%%%%

Since $\Sigma$ was assumed to have $n \leq 2c_\sigma + 1$ double points (after adding extra kinks), the complexity of $W$ is less than or equal to $2c_\sigma + 1$, completing the proof of this case. 
\smallskip

\noindent
\bit{Case 2: $\sigma^2= \pm 1$.}
\smallskip

Add at most $c_\sigma$ kinks, as in Figure \ref{kink}, to the immersion $\Sigma$ until there are algebraically zero but geometrically $n \leq 2 c_\sigma$ double points. Then the class $\lambda=\sigma + a^* \in H_2(X^\circ)$ also has square $\pm1$, and has an embedded spherical representative $\Lambda$ gotten by using the Norman trick \cite{norman} to tube the $n$ double points of a connected sum $\mathcal A^* \cs \Sigma$ over the sphere $\mathcal A$. The spheres $\Lambda$ and $\mathcal A^*$ intersect in $n$ points, whereas $\Lambda$ and $\mathcal A$ intersect geometrically once. 

As in the previous case, build a cobordism $W$ from $X$ to itself following the construction in Remark \ref{construction}, now with $\phi= R S_\lambda$ the composition of the elementary automorphisms defined in parts (b) and (c) of Definition \ref{maps}. We again claim that $W$ is the h-cobordism corresponding to the reflection $\rho_\sigma$, and so must check that conditions (1) and (2) from Remark \ref{construction} are satisfied. Note that the induced map on $H_2(X^\circ)$ sends $a^* \mapsto a$. Furthermore, setting $a,a^*=0$, the induced map sends each $\xi \in H_2(X)$ to its reflection $\xi \mp 2(\xi \cdot \sigma)\sigma$ in the class $\sigma$. Hence both conditions are satisfied. 

To bound the complexity of $W$ from above, we compute the geometric intersection number between $\mathcal A^*$ and its image $\phi(\mathcal A^*)= S_\lambda(\mathcal A)$ in $X^\circ$ (as noted in Case $1$, these are the ascending/descending spheres of the $2/3$-handle of $W$, respectively). Let $N_\Lambda$ denote a tubular neighborhood of $\Lambda$ diffeomorphic to a $2$-disk bundle over the sphere with Euler class $\pm 1$, or equivalently, a punctured copy of $\pm \mathbb{C}P^2$. The boundary $N_\Lambda \cong S^3$ is then equipped with the corresponding Hopf fibration. Let $C \subset N_\Lambda$ be a collar of $\partial N_\Lambda$. Recall from Definition \ref{maps} that the diffeomorphism $S_\lambda: X^\circ \to X^\circ$ restricts to 
\begin{enumerate}
\item complex conjugation on $N_\Lambda - C$, 
\item the identity on $X^\circ - N_\Lambda$, and
\item an isotopy $H$ on the collar $C$ from complex conjugation to the identity. 
\end{enumerate}

The neighborhood $N_\Lambda$ can be chosen small enough so that the sphere $\mathcal A$ intersects $N_\Lambda$ transversally in a single disk fiber $D$, while $\mathcal A^*$ intersects $N_\Lambda$ in the disjoint union of disk fibers $D^*_1, \dots, D^*_n$.  Let $A$ and $A^*_1, \dots, A^*_n$ denote the properly embedded annuli in which these disks intersect $C$. Identifying $C$ with $S^3 \times I$, each annulus can be thought of as a product $\delta \times I \subset S^3 \times I$ where $\delta \subset S^3$ is a circle fiber of the Hopf fibration. 

However, since complex conjugation restricted to $S^3$ sends each circle fiber to its antipodal fiber by an orientation reversing map, for each $i$ the union $H(A) \cup A^*_i$ is an immersed concordance in $S^3 \times I$ from an (oriented) Hopf link with linking number $\pm 1$\foot{Where the sign of the linking number is given by the sign of the Hopf fibration, i.e. the sign of $\sigma^2$.} to a Hopf link of opposite sign. It follows that, no matter which isotopy $H$ is chosen, the image $H(A)$ must intersect each annulus $A_i^*$ twice algebraically-- and hence at least twice geometrically. In fact, choosing the isotopy $H$ depicted in Figure \ref{iso}, this intersection number can be made \emph{exactly} two. 

Observe that $\mathcal A$ can be perturbed if necessary so that the image of the disk $D- A$ under complex conjugation on $N_\Lambda - C$ is disjoint from the union of the disks $D_i^*-A_i^*$. Therefore, the only intersection points between $\mathcal A^*$ and $S_\lambda(\mathcal A)$ are the unique point where $\mathcal A$ and $\mathcal A^*$ intersect in $X^\circ - N_\Lambda$, and the $2n$ points where the annulus $H(A)$ intersects the union of the $A_i^*$ in $C$. As $n \leq 2 c_\sigma$, this implies that the complexity of $W$ is bounded above by $4 c_\sigma +1$. 
\qed 

%%%%%%%%%% FIG 1 %%%%%%%%%%
\fig{210}{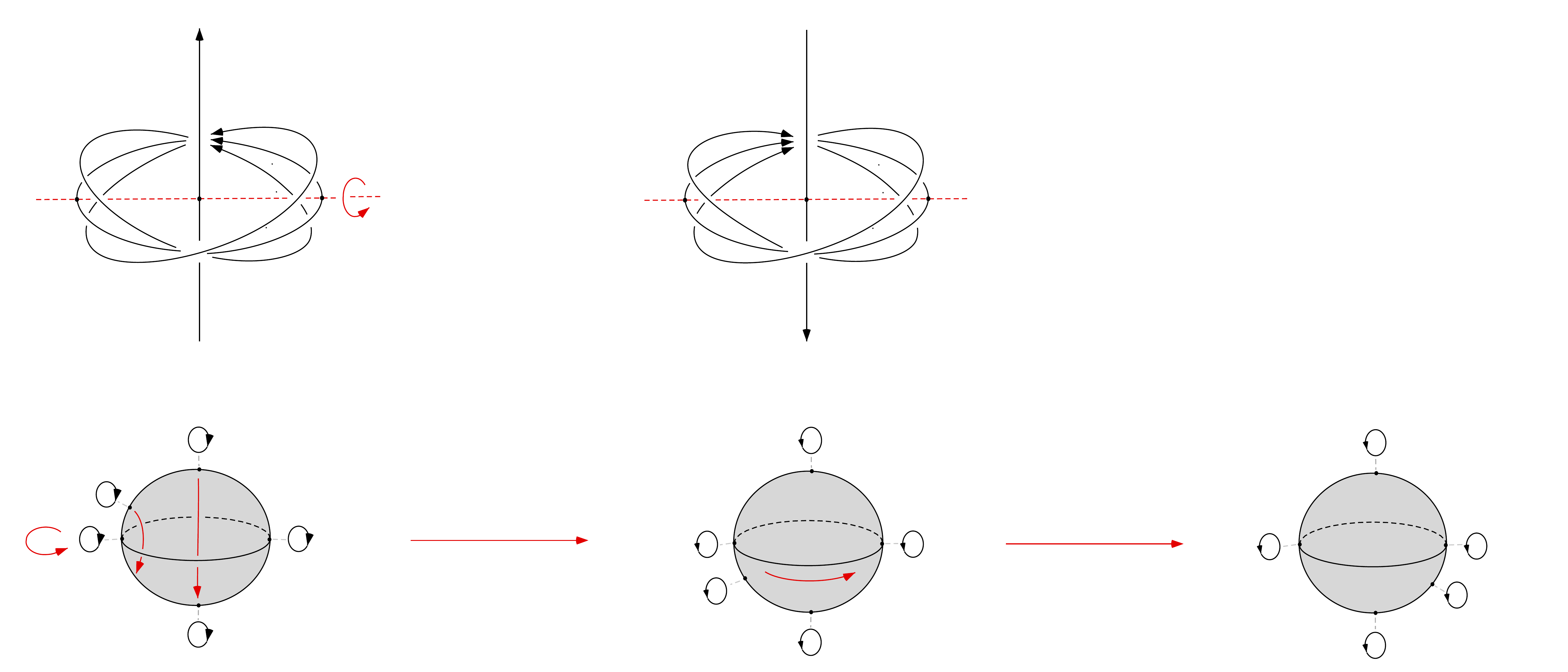}{
\put(-368,27){\small{Reflects across}}
\put(-365,16){\small{equator of $S^2$}}
\put(-150,44){\small{s}}
\put(-340,44){\small{t}}
\put(-173,27){\small{Rotates $S^2$}}
\put(-390,126){\small{$p_-$}}
\put(-390,170){\small{$p_+$}}
\put(-201,124){\small{$t(p_+)=-p_-$}}
\put(-201,169){\small{$t(p_-)=-p_+$}}
\put(-421,70){\small{$p_+$}}
\put(-421,8){\small{$p_-$}}
\caption{First, define an isotopy $t: S^3 \times I \to S^3$ which rotates by $\pi$ about the $x$-axis. At the end of the isotopy, each circle fiber is mapped to another by an orientation reversing map; however, only the fibers $p_\pm$ above the poles of $S^2$ are sent to their antipodes. This is remedied by composing with the isotopy $s: S^3 \times I \to S^3$ induced by rotating the base $S^2$ by $\pi$ about the axis that fixes the poles. Setting $H=st$, the annulus $H(A)$ intersects each $A_i^*$ exactly twice (during the isotopy $t$).}
\label{iso}}
%%%%%%%%%%%%%%%%%%%%%%%%%%%
\smallskip
Applying Theorem $1.1$ of \cite{morganszabo} in conjunction with Theorem \ref{main1} above produces some immediate results about immersed spheres in manifolds homeomorphic to rational surfaces. We establish some notation used in both of the following arguments. For each $n \in \mathbb{N}$, let $h \in H_2(\mathbb{C}P^2 \cs n \cptwobar)$ denote the class represented by $\mathbb{C}P^1 \subset \mathbb{C}P^2$ and $e_i \in H_2(\mathbb{C}P^2 \cs n \cptwobar)$ denote the class represented by the $i^{th}$ copy of $\cponebar \subset \cptwobar$. 
\smallskip

\noindent \emph{Proof of Corollary \ref{main2}.}

For each $n \in \mathbb{N}$, fix an identification $H_2(Y_n) \to H_2(\mathbb{C}P^2 \cs n \cptwobar)$ and, for brevity of notation, refer to elements in $H_2(Y_n)$ by their images. 

Assume first that $k=-1$. Then, for any $m\geq 0$ and $n=(2m+1)^2+1$, the class $y_n= (2m+1)h - \sum_{i=1}^n e_i$ is a characteristic class of square $k$ in $H_2(Y_n)$. We follow Morgan and Szab\'o's argument for Theorem $1.1$ of \cite{morganszabo} to show that the complexities $c_n$ of the isometries $\rho_n: Q_{Y_n} \to Q_{Y_n}$ reflecting in $y_n$ are unbounded. 

Since $y_n$ is characteristic of negative square, there are well-defined Seiberg-Witten invariants $SW_{Y_n, C^\pm}(y_n) \in \mathbb{Z}$ for each chamber $C^\pm \subset \{y\in H^2(Y_n; \mathbb{R}) \st y^2 =1\}$ corresponding to $y_n$ (the same two chambers also correspond to $-y_n$). By the wall-crossing formula,\foot{There are many sources giving more exposition on the Seiberg-Witten invariant, see for instance \cite{donaldson2}, \cite{LL2}, \cite{morgan}.} $$SW_{Y_n, C^+}(y_n) \pm 1=SW_{Y_n, C^-}(y_n) = \pm SW_{Y_n, C^-}(-y_n).$$ Thus $SW_{Y_n, C^+}(y_n) \not =  SW_{Y_n, \rho_n(C^+)}(\rho_n(y_n))$, since the reflection $\rho_n$ sends $y_n \mapsto -y_n$ and swaps the chambers $C^+$ and $C^-$. 

It then follows from \cite{morganszabo} that $$g(c_n) > \frac{1}{4}\big(y_n^2 -2\chi(Y_n) -3 \sigma(Y_n)\big)=n-10$$ where $g:\mathbb{Z}^+ \to \mathbb{Z}^+$ is the function given by \cite[Theorem 1.2]{morganszabo}. Since $n-10$ increases with $n$, the complexities $c_n$ must take on infinitely many values, and are thus unbounded. So by Theorem \ref{main1}, the complexity of the classes $y_n$ must grow arbitrarily large as well. A similar argument works for square $k =-2$, setting $n=(2m+1)^2+2$.  

Since the Seiberg-Witten invariant of classes with non-negative square have no chamber structure, when $k\geq 0$ we must modify our argument. Let $n=(2m+1)^2- k$ for $m$ sufficiently large so that $n>0$, and consider a new sequence of manifolds $Y'_n=Y_n \cs (k+1) \cptwobar$. The argument above shows that for the characteristic classes $y'_n= (2m+1)h - \sum_{i=1}^{n+k+1} e_i$ of square $-1$ in $H_2(Y_n')$, the complexities $c'_n$ of the isometries reflecting in $y'_n$ are unbounded. 

As $H_2(Y_n)$ sits naturally in $H_2(Y_n')$, each $y'_n$ can be thought of as the sum $y_n - \sum_{i=n+1}^{n+k+1} e_i$ where $y_n= (2m+1)h - \sum_{i=1}^n e_i$ is a characteristic class in $H_2(Y_n)$ of square $k$. Let $c_n$ denote the complexity of the reflection of $H_2(Y_n)$ in this class. Note that $c_n' \leq c_n$, since an immersed sphere in $Y_n'$ representing each class $y_n'$ can be gotten by tubing an immersed sphere in $Y_n$ representing $y_n$ to each of the $k+1$ copies of $\cponebar \subset \cptwobar$ in $Y_n'$. Thus, the complexities $c_n$ must also be unbounded. \qed 

\smallskip

\noindent \emph{Proof of Corollary \ref{main3}.}
Suppose $Y$ is homeomorphic to the rational surface $\mathbb{C}P^2 \cs n \cptwobar$, and as in the proof of Corollary \ref{main2}, fix an identification $H_2(Y) \to H_2(\mathbb{C}P^2 \cs n \cptwobar)$ so that we can refer to elements in $H_2(Y)$ by their images. Take any sequence $y_m \in H_2(Y)$ of characteristic classes with \emph{strictly increasing} square. To prove the corollary, it suffices to show that there at most finitely many terms in the sequence $y_m$ of any given complexity.  

For each $\ell_m= y_m^2+ 1$, consider the characteristic class $z_m = \sum_{i=1}^{\ell_m} e_i \in H_2(\ell_m \cptwobar)$ of square $- \ell_m$. The sum $y_m + z_m \in H_2(Y \cs \ell_m \cptwobar)$ is then also characteristic, and has square $-1$. Since the sequence $\ell_m$ increases, the argument for $k=-1$ in the proof of Corollary \ref{main2} shows that the complexity of the classes $y'_m$ must grow arbitrarily large. This forces the complexity of the classes $y_m$ to increase as well, since immersed spheres representing the $y'_m$ can be found by tubing an immersed sphere in the $Y$ summand representing $y_m$ to each $\cponebar \subset \cptwobar$ in the manifold $Y \cs \ell_m \cptwobar$. 
\qed 

\begin{remark*} Many related (and overlapping) families of homology classes of rational surfaces have been shown to have arbitrarily high complexity; Lawson's survey \cite{lawson} provides a detailed summary. Most results give lower bounds for the minimum number of \emph{positive} double points (see Fintushel and Stern \cite[Theorem 1.2]{fs:immersions}) or minimal genus (such as \cite{gqz}, \cite{L}, \cite{danny}). 

As mentioned in the introduction, Corollary \ref{main2} is implied by Ruberman \cite{danny} when each $Y_n$ is diffeomorphic to $\mathbb{C}P^2 \cs n \cptwobar$, and $k \geq 0$. Corollary \ref{main3} is also implied when $Y_n=\mathbb{C}P^2 \cs n \cptwobar$ and $n \leq 9$ by a result due to Strle \cite{saso} that there are at most finitely many \emph{reduced}\footnote{ The term ``reduced" is defined in \cite{LL}.  By \cite{LL} and \cite{L}, each class in $H_2(\mathbb{C}P^2 \cs n \cptwobar)$ of non-negative self-intersection can be sent to a reduced class by the induced map of some orientation preserving diffeomorphism of $\mathbb{C}P^2 \cs n \cptwobar$.} homology classes of any given minimal genus (and hence complexity) in a manifold $Y$ homeomorphic to a rational surface with $n \leq 9$ \cite[Proposition 14.2]{saso}.
\end{remark*}

\noindent \bit{Question.} The minimal genus of a class is bounded above by twice the complexity. However, can the difference between the complexity and the minimal genus of a class be arbitrarily large?

\end{document}